\documentclass[11pt]{amsart}
\usepackage{amssymb}
\usepackage{amscd}

\newtheorem{Thm}{Theorem}[section]
\newtheorem{Lem}[Thm]{Lemma}
\newtheorem{Cor}[Thm]{Corollary}
\newtheorem{Prop}[Thm]{Proposition}

\newcommand{\field}{k}

\begin{document}

\title[Stable endomorphism algebras of modules]
{Stable endomorphism algebras of modules over special biserial
algebras}

\author{Jan Schr\"oer}
\address{Jan Schr\"oer\newline
Department of Pure Mathematics\newline
University of Leeds\newline
Leeds LS2 9JT\newline
UK}
\email{jschroer@maths.leeds.ac.uk}

\author{Alexander Zimmermann}
\address{Alexander Zimmermann\newline
Facult\'e de Math\'ematiques et CNRS (LAMFA UMR 6140)\newline
Universit\'e de Picardie Jules Verne\newline 33, rue St
Leu\newline 80039 Amiens\newline FRANCE}
\email{alexander.zimmermann@u-picardie.fr}

\thanks{
We gratefully acknowledge support from the Volkswagen
Foundation (RIP Program at Oberwolfach).}

\begin{abstract}
We prove that the stable endomorphism algebra of a module without
self-extensions over a special biserial algebra is a gentle
algebra. In particular, it is again special biserial. As a
consequence, any algebra which is derived equivalent to a gentle
algebra is gentle.
\end{abstract}

\maketitle

\begin{center}
Dedicated to Idun Reiten on the occasion of her 60th birthday
\end{center}




\section{Introduction}

Let $\field$ be an algebraically closed field, and let $A$ be a 
$\field$-algebra. 
One of the most important classes of tame algebras is 
the class of special biserial algebras. 
These algebras occur naturally in many different contexts, 
see the introduction of \cite{Sc3} or \cite{Sc4} for more details. 
We shall prove that for a special biserial algebra $A$ and 
any $A$-module $M$ without self-extensions the endomorphism algebra 
of $M$ modulo the ideal of 
endomorphisms factoring through a projective 
$A$-module is not only special biserial again, but is even 
a gentle algebra.
Gentle algebras are special biserial
algebras satisfying certain minimality conditions.

The result is interesting and surprising in its own right, since 
for arbitrary
algebras this quotient of the endomorphism algebra of modules without 
self-extensions will be very different from the original algebra. 
Nevertheless, our main motivation came from different considerations, namely 
the theory of derived equivalences. 
Let ${\rm mod}(A)$ be the category of finite-dimensional right 
$A$-modules and let
$D^b(A)$ be the derived category of bounded complexes of
finitely generated $A$-modules.
Two algebras $A$ and $B$ are said to be derived equivalent if
the categories $D^b(A)$ and $D^b(B)$ are equivalent as triangulated
categories. The stable module category $\underline{\rm mod}(A)$ is 
the category with the same
objects as ${\rm mod}(A)$, and morphisms are equivalence classes of 
$A$-module 
homomorphisms modulo those which factor through a projective module.  

Rickard shows that two algebras $A$ and $B$ are derived equivalent if
and only if there exists a so-called tilting complex $T$ in $D^b(A)$ 
such that
$B \simeq {\rm End}_{D^b(A)}(T)$, see \cite{derbuch},\cite{R}.
Now, Happel proved  in \cite{H} that there is a full and faithful embedding
$F: D^b(A) \longrightarrow \underline{\rm mod}(RA)$
of triangulated categories, where $RA$ is the repetitive algebra of $A$.
Since for a tilting complex $T$ one gets ${\rm Hom}_{D^b(A)}(T,T[1])=0$, 
the algebra $B$ is isomorphic to the endomorphism algebra in 
$\underline{\rm mod}(RA)$ of a module $M=FT$ without self-extensions.
Thus, studying stable endomorphism algebras of modules 
without self-extensions is of central importance 
when one is interested in derived equivalences.  

In our case, since an algebra is gentle if and only if its 
repetitive algebra is special 
biserial, (see \cite{PS} and also \cite{Ri2} and \cite{Sc1}), we get 
as a corollary to our main result that the class of gentle
algebras is closed under derived equivalences.
At the moment, examples of this kind are very rare and clearly show that
gentle algebras deserve much attention.

\medskip

To be more precise, let us recall some definitions:
Let $Q$ be a quiver.
For an arrow $\alpha: a \to b$, let $s(\alpha) = a$ be its starting point
and $e(\alpha) = b$ its end point.
A {\it path} of length $n \geq 1$ in $Q$ is a sequence
$\alpha_1 \cdots \alpha_n$
of arrows such that $e(\alpha_i) = s(\alpha_{i+1})$ for all
$1 \leq i \leq n-1$.
A {\it relation} for $Q$ is a non-zero $\field$-linear combination of
paths of length at least 2 having the same starting point and the same
end point.
Let $\rho$ be a set of relations for $Q$.
Then $(Q,\rho)$ is called {\it special biserial} if
the following hold:
\begin{enumerate}

\item[(1)] Any vertex in $Q$ is the starting point of at most two arrows
and also the end point of at most two arrows;

\item[(2)] Given an arrow $\beta$, there is at most one arrow $\alpha$ with
$e(\alpha) = s(\beta)$ and $\alpha\beta \notin \rho$,
and there is at most one arrow $\gamma$ with $e(\beta) = s(\gamma)$ and
$\beta\gamma \notin \rho$;

\item[(3)] Each infinite path in $Q$ contains a subpath which is
in $\rho$.
\end{enumerate}
A special biserial pair $(Q,\rho)$ is {\it gentle} if additionally the
following hold:
\begin{enumerate}

\item[(4)] All elements in $\rho$ are paths of length 2;

\item[(5)] Given an arrow $\beta$, there is at most one arrow $\alpha'$
with $e(\alpha') = s(\beta)$ and
$\alpha'\beta \in \rho$,
and there is at most one arrow $\gamma'$ with $e(\beta) = s(\gamma')$ and
$\beta\gamma' \in \rho$.

\end{enumerate}
A $\field$-algebra is called {\it special biserial}, or {\it gentle},
if it is Morita equivalent to an algebra $\field Q/(\rho)$
for $(Q,\rho)$ special biserial, or gentle, respectively.
In these cases, $\field Q/(\rho)$ is finite-dimensional if and only if
$Q$ contains only finitely many vertices.
Here $\field Q$ is the path algebra of $Q$, and $(\rho)$ is the ideal
generated by the elements in $\rho$.
Recall that $\field Q$ has as a basis the set of all paths in
$Q$ including a path $e_i$ of length 0 for each vertex $i$ in $Q$.
By `modules' we always mean finitely generated right modules.

Here is our main result.

\begin{Thm}\label{main}
Let $A$ be a special biserial algebra, and let $M$ be an $A$-module. If
${\rm Ext}^1_A(M,M)=0$, then $\underline{\rm End}_A(M)$ is a gentle
algebra.
\end{Thm}

Let $T[i]$ be the usual shift of a complex $T$ in $D^b(A)$ by
$i$ degrees to the left (we adopt the convention of \cite{derbuch}).
Then, this theorem has the following consequence.

\begin{Cor}\label{maincor}
Let $A$ be a finite-dimensional gentle algebra, and let $T$ be a complex in
$D^b(A)$.
If ${\rm Hom}_{D^b(A)}(T,T[1])=0$, then ${\rm End}_{D^b(A)}(T)$ is
a gentle algebra.
In particular, any algebra $B$ which is derived equivalent to $A$
is gentle, and (up to Morita equivalence) there are only finitely many
such algebras $B$.
\end{Cor}

Classical examples of gentle algebras are hereditary algebras of type
$\mathbb{A}_n$ and $\tilde{\mathbb{A}}_n$.
It was known before that all algebras which are derived equivalent
to these examples are gentle again.
For a complete classification of the derived equivalence classes
of these examples see \cite{AH} and \cite{AS}.
Recently, Vossieck classified all algebras $A$ such that $D^b(A)$ is
discrete.
In this case, it turns out that $A$ is either derived hereditary of
Dynkin type or a gentle algebra, see \cite{V}.
Using Vossieck's result, the derived equivalence classes of
algebras with discrete derived category have been classified in \cite{BGS}.

The paper is organised as follows.
In Section \ref{basics} we recall definitions and give a survey of
known results.
The main theorem and its corollary are proved in Section \ref{proofs}.
Finally, in Section \ref{examples} we give some examples.

Although we often write maps on the left hand side,
we compose them as if they were on the right.
Thus the composition of a map $\theta$
followed by a map $\phi$ is denoted $\theta\phi$.




\section{Known facts on special biserial algebras}\label{basics}

In this section, we recall some basic facts on special biserial algebras,
string modules and homomorphisms between string modules.
As a main reference we use \cite{BR}, see also there for further references.

Let $Q = (Q_0,Q_1)$ be a quiver with set of vertices $Q_0$ and set of
arrows $Q_1$.
Let $\rho$ be a set of relations for $Q$
such that $(Q,\rho)$ is special biserial, and let $A = \field Q/(\rho)$
be the corresponding special biserial algebra.
Without loss of generality we can assume that $\rho$ contains only
zero relations, i.e. paths, and relations
of the form $p-q$ with $p$ and $q$ paths, which are called
{\it commutativity relations.}
If $r = p-q$ is a commutativity relation, then we say that $p$ and $q$
are contained in $r$.
By $\rho^+$ we denote the set $\rho$ together with all paths which
are contained in commutativity relations in $\rho$.

Given an arrow $\alpha$ in $Q$
denote by $\alpha^-$ a formal inverse where $s(\alpha^-) = e(\alpha)$ and
$e(\alpha^-) = s(\alpha)$.
Also let $(\alpha^-)^- = \alpha$.
The set of formal inverses of arrows is denoted by $Q_1^-$.
We use the word `arrow' only for element in $Q_1$, we refer to elements
in $Q_1^-$ always as inverse arrows.
A {\it string} (for $(Q,\rho)$) of length $n \geq 1$ is a sequence
$c_1 \cdots c_n$ of arrows or inverse arrows with the following properties:
\begin{enumerate}

\item[(1)] $e(c_i) = s(c_{i+1})$ for $1 \leq i \leq n-1$;

\item[(2)] $c_i \not= c_{i+1}^-$ for $1 \leq i \leq n-1$;

\item[(3)] For all $1 \leq i < j \leq n$ neither $c_i c_{i+1} \cdots c_j$
nor $c_j^- \cdots c_{i+1}^- c_i^-$ are contained in $\rho^+$.

\end{enumerate}
Let $C = c_1 \cdots c_n$ be a string.
Define $s(C) = s(c_1)$ and $e(C) = e(c_n)$.
Let $C^- = c_n^- \cdots c_1^-$ be the inverse string of $C$.
Additionally, we define for every vertex $i$ in $Q$ two strings
$1_{(i,1)}$ and $1_{(i,-1)}$ of length 0 with
$s(1_{(i,t)}) = e(1_{(i,t)}) = i$ and $1_{(i,t)}^- = 1_{(i,-t)}$
for $t = -1,1$.
The length of an arbitrary string $C$ is denoted by $|C|$.
We call $C$ {\it direct} if $|C| = 0$ or
$C = c_1 \cdots c_n$ with $c_i \in Q_1$ for all $i$, and $C$ is called
{\it inverse} if $C^-$ is direct.

Similarly as in \cite{BR}, one can define two maps
$\sigma,\epsilon: Q_1 \to \{ -1,1 \}$ which satisfy the following
properties:
\begin{enumerate}

\item[(1)] If $\alpha_1 \not= \alpha_2$ are arrows with
$e(\alpha_1) = e(\alpha_2)$, then $\epsilon(\alpha_1) =
- \epsilon(\alpha_2)$;

\item[(2)] If $\beta_1 \not= \beta_2$ are arrows with
$s(\beta_1) = s(\beta_2)$, then $\sigma(\beta_1) = - \sigma(\beta_2)$;

\item[(3)] If $\alpha$ and $\beta$ are arrows with $e(\alpha) = s(\beta)$
and $\alpha\beta \notin \rho^+$, then
$\epsilon(\alpha) = - \sigma(\beta)$;

\item[(4)] Let $i$ be a vertex such that there is only one arrow
$\alpha$ with $e(\alpha) = i$, and only one arrow $\beta$ with
$s(\beta) = i$.
If $\alpha\beta \in \rho^+$, then $\epsilon(\alpha) = \sigma(\beta)$.

\end{enumerate}
Note that $\sigma$ and $\epsilon$ are defined in \cite{BR} without
condition (4).
The set of all strings is denoted by ${\mathcal S}$.
We extend $\sigma$ and $\epsilon$ to maps
\[
\sigma,\epsilon: {\mathcal S} \to \{ -1,1\}
\]
as follows:
For an arrow $\alpha$ define $\sigma(\alpha^-) = \epsilon(\alpha)$ and
$\epsilon(\alpha^-) = \sigma(\alpha)$.
If $C = c_1 \cdots c_n$ is a string of length $n \geq 1$, then let
$\sigma(C) = \sigma(c_1)$ and $\epsilon(C) = \epsilon(c_n)$.
For a string $1_{(i,t)}$ of length 0 define $\sigma(1_{(i,t)}) = t$
and $\epsilon(1_{(i,t)}) = -t$.

Let $C = c_1 \cdots c_n$ and $D = d_1 \cdots d_m$ be strings of length
at least 1.
If $CD = c_1 \cdots c_n d_1 \cdots d_m$ is a string, then we say that the
concatenation of $C$ and $D$ is defined.
For an arbitrary string $C$ let $1_{(s(C),t)}C = C$ if
$\sigma(C) = t$, and let $C1_{(e(C),t)} = C$ if $\epsilon(C) = -t$.
Otherwise, the concatenation with a string of length 0 is not defined.

One uses $\sigma$ and $\epsilon$ to give a certain `orientation' to
strings.
One of the main properties of $\sigma$ and $\epsilon$ is the following:

\begin{Lem}
If $C$ and $D$ are strings such that the concatenation $CD$ is defined,
i.e. CD is again a string, then $e(C) = s(D)$ and
$\epsilon(C) = - \sigma(D)$.
\end{Lem}

For each string $C$, we define a map
\[
v_C: \{ 1, \cdots, |C|+1 \} \longrightarrow Q_0
\]
as follows.
If $C = 1_{(i,t)}$ is a string of length 0, then let
$v_C(1) = i$.
If $C = c_1 \cdots c_n$ is of length at least 1, then let
\[
v_C(i) = \left\{ \begin{array}{r@{\quad : \quad}l}
s(c_i) & 1 \leq i \leq |C|; \\
e(c_n) & i = |C|+1.
\end{array} \right.
\]
For each string $C$ we
construct an $A$-module $M(C)$ as follows.
First, assume that $C = c_1 \cdots c_n$ is a string of length $n \geq 1$.
Fix a basis $\{ z_1, \cdots ,z_{n+1} \}$ of $M(C)$.
Let $e_j$ be a path of length 0 in $Q$.
Then define
\[
z_i \cdot e_j = \left\{ \begin{array}{r@{\quad : \quad}l}
z_i & i \in v_C^{-1}(j); \\
0 & $otherwise.$
\end{array} \right.
\]
Given an arrow $\alpha$ in $Q$ let
\[
z_i \cdot \alpha = \left\{ \begin{array}{r@{\quad : \quad}l}
z_{i+1} & \alpha=c_i, 1 \leq i \leq n;\\
z_{i-1} & \alpha^- = c_{i-1}, 2 \leq i \leq n+1;\\
0 & $otherwise.$
\end{array} \right.
\]
Next, assume $C = 1_{(i,t)}$ is a string of length 0.
Then let $\{ z_1 \}$ be a basis of $M(C)$.
For each path $e_j$ of length 0 in $Q$ let
\[
z_1 \cdot e_j = \left\{ \begin{array}{r@{\quad : \quad}l}
z_1 & j=i; \\
0 & $otherwise.$
\end{array} \right.
\]
For all $\alpha$ in $Q_1$ define $z_1 \cdot \alpha = 0$.

Modules of the form $M(C)$ are called {\it string modules}, and
$\{ z_1, \cdots, z_{n+1} \}$ is
called the {\it canonical basis} of $M(C)$.
The construction of string modules goes back to Gelfand and Ponomarev, see
\cite{GP} and also \cite{BR}.
One can easily check that string modules are always indecomposable.
For strings $C_1$ and $C_2$ define $C_1 \sim C_2$ if $C_1 = C_2$ or
$C_1 = C_2^-$.
Two string modules $M(C_1)$ and $M(C_2)$ are isomorphic if and only
if $C_1 \sim C_2$.

For each vertex $i$, let $S_i = M(1_{(i,1)})$.
Then $\{ S_i \mid i \in Q_0 \}$ is a complete set of isomorphism
classes of simple $A$-modules.
Furthermore, let $P_i$ be the corresponding indecomposable projective
$A$-module with ${\rm top}(P_i) = S_i$.

For a string $C$ define
${\mathcal P}(C) = \{ (D,E,F) \mid D,E,F \in {\mathcal S} \;{\rm and}\;
DEF = C \}$.
We call $(D,E,F)$ in ${\mathcal P}(C)$ a {\it factor string} of $C$ if the
following hold:
\begin{enumerate}

\item[(1)] Either $|D|= 0$ or $D = d_1 \cdots d_n$ where
$d_n \in Q_1^-$;

\item[(2)] Either $|F|= 0$ or $F = f_1 \cdots f_m$ where
$f_1 \in Q_1$.

\end{enumerate}
Dually, we call $(D,E,F)$ a {\it substring} of $C$ if the
following hold:
\begin{enumerate}

\item[(1)] Either $|D|= 0$ or $D = d_1 \cdots d_n$ where
$d_n \in Q_1$;

\item[(2)] Either $|F|= 0$ or $F = f_1 \cdots f_m$ where
$f_1 \in Q_1^-$.

\end{enumerate}
Let ${\rm fac}(C)$ be the set of factor strings of $C$,
and by ${\rm sub}(C)$ we denote the set of substrings of $C$.
If $C_1$ and $C_2$ are strings, then we
call a pair $((D_1,E_1,F_1),(D_2,E_2,F_2))$ in
${\rm fac}(C_1) \times {\rm sub}(C_2)$ {\it admissible} if $E_1 \sim E_2$.
Denote the set of admissible pairs in ${\rm fac}(C_1) \times {\rm sub}(C_2)$
by ${\mathcal A}(C_1,C_2)$.

For each $a = ((D_1,E_1,F_1),(D_2,E_2,F_2))$ in ${\mathcal A}(C_1,C_2)$
we define a homomorphism
\[
f_a : M(C_1) \longrightarrow M(C_2)
\]
as follows:
If $E_1 = E_2$, then let
\[
f_a(z_{|D_1|+i}) = z_{|D_2|+i}
\]
for $1 \leq i \leq |E_1|+1$, and
all other canonical basis vectors of $M(C_1)$ are mapped to 0.
If $E_1 = E_2^-$, then let
\[
f_a(z_{|D_1|+i}) = z_{|D_2E_2|+2-i}
\]
for $1 \leq i \leq |E_1|+1$, and all other canonical basis vectors are
mapped to 0.
Such homomorphisms are called {\it graph maps}.
It is proved in \cite{CB} that the graph maps
$\{ f_a \mid a \in {\mathcal A}(C_1,C_2) \}$ form a basis
of the homomorphism space ${\rm Hom}_A(M(C_1),M(C_2))$.

There is the following multiplicative behaviour of graph maps:
The composition $f_a f_b$ of graph maps is either 0 or a graph map.
This fact is very important and has numerous applications.

If $a = ((D_1,E_1,F_1),(D_2,E_2,F_2))$ is admissible, then we call
$a$ and $f_a$ {\it oriented} if $E_1 = E_2$. We say that $a$ and
$f_a$ are {\it left-sided}, or {\it right-sided} if they are
oriented and $|D_1| = |D_2| = 0$, or $|F_1| = |F_2| = 0$,
respectively. In case $a$ is left-sided or right-sided, one calls
$a$ and $f_a$ {\it one-sided}. Finally, $a$ and $f_a$ are {\it
weakly one-sided} if $a$ or $((D_1,E_1,F_1),(F_2^-,E_2^-,D_2^-))$
are one-sided, and they are {\it two-sided} if they are not
weakly one-sided.

Let $C_1 = D_1E_1F_1$ and $C_2 = D_2E_2F_2$.
A graph map $f_a$ can always be transformed to an oriented
graph map.
This is done by composing $f_a$ with one of the canonical
isomorphisms $M(C_1^-) \to M(C_1)$ or $M(C_2) \to M(C_2^-)$.
Note that these are both graph maps.
Since there are two possibilities, we define
$a(l) = a(r) = a$ if $a$ is oriented, and
$a(l) = ((F_1^-,E_1^-,D_1^-),(D_2,E_2,F_2))$ and
$a(r) = ((D_1,E_1,F_1),(F_2^-,E_2^-,D_2^-))$, otherwise.
It follows that $a(l)$ is one-sided if and only if $a(r)$
is one-sided if and only if $a$ is weakly one-sided.
The next lemma follows directly from the definitions.

\begin{Lem}\label{multiplicative}
Let $f_a$ and $f_b$ be weakly one-sided graph maps.
If $f_{a(l)}$ and $f_{b(r)}$ are both left-sided, or both right-sided,
then $f_{a(l)}f_{b(r)}$ is again left-sided, or right-sided, respectively.
In particular, $f_af_b$ is again weakly one-sided.
Otherwise, $f_af_b$ is 0 or two-sided.
\end{Lem}

In \cite{Sc2} we use the term `perfect' instead of `oriented'.
Note also that we use left modules in \cite{Sc2} instead of right modules.
Thus strings are written in the opposite way.
One of the important steps in our work is to prove a
stable version of the following proposition.

\begin{Prop}\label{oldversion}
Let $A = \field Q/(\rho)$ be a special biserial algebra, and let
$C_i$, $1 \leq i \leq t$, be strings for $(Q,\rho)$.
If each $a \in {\mathcal A}(C_i,C_j)$, $1 \leq i,j \leq t$,
is weakly one-sided, then ${\rm End}_A(\bigoplus_{i=1}^t M(C_i))$
is a gentle algebra.
\end{Prop}

\begin{proof}
For $A$ a string algebra, this statement is proved in
\cite[Proposition 4.8]{Sc2}.
Now assume $A$ is special biserial, and let $p-q$ be a commutativity
relation in $\rho$.
If $M$ is a direct sum of string modules, then $M \cdot p$ and
$M \cdot q$ are both 0.
Thus we can define $B = A/(p_i, i \in I)$ where the $p_i$ are all paths
which are contained in commutativity relations.
Clearly, $B$ is a string algebra, and ${\rm End}_A(M)$ is isomorphic
to ${\rm End}_B(M)$.
\end{proof}




\section{Proof of the main result}\label{proofs}

Throughout this section, we assume that $A = \field Q/(\rho)$
is special biserial, with $\rho$ containing only zero relations or
commutativity relations.

By $\underline{\rm Hom}_A(M,N)$ we denote the $A$-module homomorphisms
from $M$ to $N$ modulo the homomorphisms which factor through
projective $A$-modules.
Let $\underline{\rm mod}(A)$ be the stable category of finite-dimensional
$A$-modules.
Here the objects are the same as in ${\rm mod}(A)$, but the morphism
space from $M$ to $N$ is $\underline{\rm Hom}_A(M,N)$.
For a homomorphism $f$ let $\underline{f}$ be the corresponding
morphism in $\underline{\rm mod}(A)$.

It is known that any finite-dimensional
indecomposable $A$-module is a string module, a
band module, or a non-serial projective-injective module, see \cite{BR}.
If $P$ is an indecomposable non-serial projective-injective
$A$-module, then the radical
${\rm rad}(P)$ is isomorphic to a string module $M(CD)$ with $C$ a direct
and $D$ an inverse string.
Similarly, the socle factor $P/{\rm soc}(P)$ is isomorphic to a string module
$M(C'D')$ with $C'$ an inverse and $D'$ a direct string.

If $N$ is a band module, then ${\rm Ext}^1_A(N,N) \not= 0$ since
band modules lie in homogeneous tubes of the Auslander-Reiten quiver
of $A$, i.e. $\tau_A N = N$, where $\tau_A$ is the Auslander-Reiten
translation.
For basic facts on Auslander-Reiten theory we refer to
\cite{ARS} or \cite{Ri1}.

So, if $M$ is an $A$-module with ${\rm Ext}_A^1(M,M) = 0$, then
$M$ is a direct sum of string modules and of non-serial
projective-injectives.
Since we are interested in stable homomorphisms, it is
enough to consider only string modules.

For the rest of this section, let
$C_1$ and $C_2$ be (not necessarily different) strings for $(Q,\rho)$.

\begin{Lem}\label{facprojinj}
Let $f: M(C_1) \to M(C_2)$ be a homomorphism which factors through
a direct sum of non-serial projective-injective $A$-modules.
Then $f$ is a linear combination of two-sided graph maps.
\end{Lem}

\begin{proof}
Clearly, it is enough to prove the lemma for the case that $f$ factorizes
through an indecomposable non-serial projective-injective
$A$-module.
Thus, let $P$ be an indecomposable non-serial projective-injective
$A$-module, and let $f_1: M(C_1) \to P$ and $f_2: P \to M(C_2)$ homomorphisms
such that $f = f_1f_2$.
By $\iota: {\rm rad}(P) \to P$ we denote the canonical radical inclusion,
and by
$\pi: P \to P/{\rm soc}(P)$ the projection of $P$ onto its socle factor.
Clearly, $f_1$ must factor through ${\rm rad}(P)$, and $f_2$ must factor
through $P/{\rm soc}(P)$.
Thus there are homomorphisms $g$ and $h$ such that $f_1 = g\iota$ and
$f_2 = \pi h$.
The following picture describes the situation:

\unitlength1.0cm
\begin{picture}(11,5.4)

\put(0.6,4.5){$M(C_1)$}
\put(8.6,4.5){$M(C_2)$}
\put(2.2,4.6){\vector(1,0){6}}
\put(5,4.8){$f$}

\put(1.1,4.1){\vector(0,-1){0.9}}
\put(1.2,3.7){$g$}

\put(9.1,3.4){\vector(0,1){0.7}}
\put(9.2,3.6){$h$}

\put(0.5,2.3){\circle*{0.1}}
\put(0.5,2.2){\vector(0,-1){0.8}}\put(0,1.7){$R_1$}
\put(1.7,2.3){\circle*{0.1}}
\put(1.7,2.2){\vector(0,-1){0.8}}\put(1.8,1.7){$R_2$}
\put(0.5,1.3){\circle*{0.1}}
\put(0.6,1.2){\vector(1,-1){0.4}}\put(0.5,0.8){$\alpha$}
\put(1.7,1.3){\circle*{0.1}}
\put(1.6,1.2){\vector(-1,-1){0.4}}\put(1.5,0.8){$\beta$}
\put(1.1,0.7){\circle*{0.1}}

\put(3.1,2.1){$\iota$}
\put(2.8,2){\vector(1,0){0.7}}

\put(5.1,2.9){\circle*{0.1}}
\put(5,2.8){\vector(-1,-1){0.4}}\put(4.4,2.6){$\gamma$}
\put(5.2,2.8){\vector(1,-1){0.4}}\put(5.5,2.6){$\delta$}
\put(4.5,2.3){\circle*{0.1}}
\put(4.5,2.2){\vector(0,-1){0.8}}\put(4,1.7){$R_1$}
\put(5.7,2.3){\circle*{0.1}}
\put(5.7,2.2){\vector(0,-1){0.8}}\put(5.8,1.7){$R_2$}
\put(4.5,1.3){\circle*{0.1}}
\put(4.6,1.2){\vector(1,-1){0.4}}\put(4.5,0.8){$\alpha$}
\put(5.7,1.3){\circle*{0.1}}
\put(5.6,1.2){\vector(-1,-1){0.4}}\put(5.5,0.8){$\beta$}
\put(5.1,0.7){\circle*{0.1}}

\put(7.1,1.9){$\pi$}
\put(6.8,1.8){\vector(1,0){0.9}}

\put(9.1,2.9){\circle*{0.1}}
\put(9,2.8){\vector(-1,-1){0.4}}\put(8.4,2.6){$\gamma$}
\put(9.2,2.8){\vector(1,-1){0.4}}\put(9.5,2.6){$\delta$}
\put(8.5,2.3){\circle*{0.1}}
\put(8.5,2.2){\vector(0,-1){0.8}}\put(8,1.7){$R_1$}
\put(9.7,2.3){\circle*{0.1}}
\put(9.7,2.2){\vector(0,-1){0.8}}\put(9.8,1.7){$R_2$}
\put(8.5,1.3){\circle*{0.1}}
\put(9.7,1.3){\circle*{0.1}}

\put(0.6,0){${\rm rad}(P)$}
\put(4.9,0){$P$}
\put(8.4,0){$P/{\rm soc}(P)$}

\end{picture}

\vspace{0.5cm}

Now $g = \sum_i \lambda_i f_{a_i}$ and $h = \sum_j \mu_j f_{b_j}$
where $a_i$ is of the form
\[
((D_{1i},E_{1i}{E_{1i}'}^-,F_{1i}),(D_{2i},E_{2i},F_{2i}))
\]
with $E_{1i}$ and $E_{1i}'$ being direct strings,
and $b_j$ is of the form
\[
((D_{3j},E_{3j}^-E_{3j}',F_{3j}),(D_{4j},E_{4j},F_{4j}))
\]
with $E_{3j}$ and $E_{3j}'$ being direct strings.
Thus we get
\[
f = \sum_{i,j} \lambda_i\mu_j f_{a_i} \iota \pi f_{b_j}.
\]
It is straightforward to check that for all $i,j$ the composition
$f_{a_i} \iota \pi f_{b_j}$
is either 0 or a sum of either one or two graph maps which are two-sided.
\end{proof}

\begin{Lem}\label{basis}
Let $f_{a_i}: M(C_1) \to M(C_2)$, $1 \leq i \leq s$,
be pairwise different graph maps,
which are weakly one-sided, with $\underline{f_{a_i}} \not= 0$.
Then the $\underline{f_{a_i}}$ are linear independent in
$\underline{\rm Hom}_A(M(C_1),M(C_2))$.
\end{Lem}

\begin{proof}
Assume that $\sum_{i=1}^s \lambda_i \underline{f_{a_i}} = 0$
with $\lambda_i \not= 0$ for some $i$, and define
$f = \sum_{i=1}^s \lambda_i f_{a_i}$.
Thus there is a projective $A$-module $P = \bigoplus_{i=1}^l P(i)$,
$P(i)$ indecomposable for all $i$, and homomorphisms
$f_1 = [f_{11}, \cdots, f_{1l}]: M(C_1) \to P$ and
$f_2 = [f_{21}, \cdots, f_{2l}]^t: P \to M(C_2)$ such that $f = f_1f_2$.
Thus $f_1f_2 = \sum_{i=1}^l f_{1i}f_{2i}$.
If $P(i)$ is not a non-serial projective-injective for some $i$,
then $f_{1i}$ and $f_{2i}$ are both linear combinations of graph maps.
By the multiplicativity of graph maps, we get that $f_{1i}f_{2i}$
is a linear combination of graph maps as well, and each of these factors
through $P(i)$.
Otherwise, if $P(j)$ is a non-serial projective-injective for some j,
Lemma \ref{facprojinj} yields that $f_{1j}f_{2j}$ is a linear
combination of two-sided graph maps.
But different graph maps are linear independent, and the
$f_{a_i}$ are by assumption weakly one-sided.
So all $f_{a_i}$ with $\lambda_i \not= 0$ must factor through a
projective, a contradiction.
\end{proof}

\begin{Lem}\label{stablemultiplicative}
Let $f_a$ and $f_b$ be weakly one-sided graph maps.
If $f_{a(l)}$ and $f_{b(r)}$ are both left-sided, or both right-sided, and
if $\underline{f_a} \not= 0 \not= \underline{f_b}$, then
$\underline{f_af_b} \not= 0$.
\end{Lem}

\begin{proof}
Without loss of generality assume that $f_a$ and $f_b$ are both
right-sided.
Define $f_c = f_af_b$.
For the sake of brevity we just write `1' instead of
$1_{(e(E),-\epsilon(E))}$.
Thus $c$ is of the form
\[
((D_1,E,1),(D_2,E,1)).
\]
Assume that $\underline{f_a} \not= 0 \not= \underline{f_b}$.
Note that this implies that $M(D_1E)$ and $M(D_2E)$ are both not
projective.
We have to show that $\underline{f_af_b} \not= 0$.

Assume that $f_c = f_1f_2$ is a factorization of $f$ through a
projective module $P$.
Since $f_c$ is weakly one-sided, and since different graph maps are
linear independent, it follows from Lemma \ref{facprojinj}
that we can assume that $P$ does not contain non-serial projective-injective
direct summands.
Using the multiplicativity of graph maps, we can assume without
loss of generality that $P = M(L^-RE)$ with $L$ and
$R$ direct strings, $|R| \geq 1$, and
$f_1$ and $f_2$ right-sided graph maps.
In particular, $E$ must be a direct string.
Note that $D_2 = D_2'R$ for some string $D_2'$.
We know that $a$ is of the form
\[
((D_{11},D_{12}E,1),(D_3,D_{12}E,1)),
\]
and $b$ is of the form
\[
((D_4,D_{22}E,1),(D_{21},D_{22}E,1)).
\]
Now let $D_{122}$ be a direct string of maximal length such that
$D_{12} = D_{121}D_{122}$ for some string $D_{121}$, and let
$D_{222}$ be a maximal direct string such that
$D_{22} = D_{221}D_{222}$ for some string $D_{221}$.
Since $f_c = f_af_b$, we get that $|D_{12}| = 0$ or $|D_{22}| = 0$.

If $|D_{22}| = 0$, then $f_b$ clearly factors through
$P = M(L^-RE)$, a contradiction.
Thus, assume $|D_{22}| \geq 1$.
This implies $|D_{12}| = 0$.
If now $|R| \leq |D_{222}|$, then $f_a$ factors through $P$.
Otherwise, if $|R| > |D_{222}|$, then $D_{22} = D_{222}$, since
$D_2 = D_2'R$ and since $b$ is admissible.
This implies that $f_b$ factors through $P$.
This finishes the proof.
\end{proof}

\begin{Prop}\label{newversion}
Let $A = \field Q/(\rho)$ be a special biserial algebra, and let
$C_i$, $1 \leq i \leq t$, be strings for $(Q,\rho)$.
If each $f_a$, $a \in {\mathcal A}(C_i,C_j)$, $1 \leq i,j \leq t$,
is weakly one-sided or factors through a projective $A$-module,
then $\underline{\rm End}_A(\bigoplus_{i=1}^t M(C_i))$
is a gentle algebra.
\end{Prop}

\begin{proof}
Since we are interested in algebras only up to Morita equivalence,
we can assume that the $M(C_i)$ are pairwise not isomorphic.
We can also assume that the $M(C_i)$ are not projective.
For the sake of brevity let $M = \bigoplus_{i=1}^t M(C_i)$.
By definition $\underline{\rm End}_A(M) = {\rm End}_A(M)/{\mathcal P}$
where ${\mathcal P}$ is the ideal of all endomorphisms which factor
through projectives.

The graph maps in ${\rm Hom}_A(M(C_i),M(C_j))$, $1 \leq i,j \leq t$, form a
basis of ${\rm End}_A(M)$.
By our assumption, we know that all two-sided graph maps are
contained in ${\mathcal P}$.
This implies that the weakly one-sided graph maps generate
$\underline{\rm End}_A(M)$ as a vectorspace.
Let ${\mathcal B}$ be the set of weakly one-sided graph maps $f_a$ in
${\rm End}_A(M)$ which satisfy $\underline{f_a} \not= 0$.
It follows form Lemma \ref{basis} that
$\underline{\mathcal B} = \{ \underline{f_a} \mid f_a \in {\mathcal B} \}$
is a basis of $\underline{\rm End}_A(M)$.
This basis behaves again multiplicative.
Namely, let $f_a$ and $f_b$ be in ${\mathcal B}$.
Then it follows from Lemma \ref{multiplicative} that
$\underline{f_af_b} = 0$,
or  $\underline{f_af_b} \in \underline{\mathcal B}$.

Now we can compute the quiver with relations of the algebra
$\underline{\rm End}_A(M)$.
The vertices of this quiver are the identity maps $M(C_i) \to M(C_i)$.
The arrows are the maps $\underline{f_a}$ in $\underline{\mathcal B}$ which
are not of the form $\underline{f_b}\underline{f_c}$ with
$\underline{f_b}$ and $\underline{f_c}$ non-invertible elements in
$\underline{\mathcal B}$.
Now we proceed as in the proof of \cite[Proposition 4.8]{Sc2}.
As in \cite[Lemmas 4.1, 4.2 and Corollary 4.3]{Sc2} we get that there are
at most two arrows ending and at most two arrows starting at each vertex
of the quiver of $\underline{\rm End}_A(M)$.
Thus we get property (1) in the definition of a gentle algebra.
Then property (2) and (5) are proved as in \cite[Lemma 4.4, 4.5]{Sc2}.
However, to get (5) we have to use also Lemma \ref{stablemultiplicative}.
Since $M$ is a finite-dimensional module, also $\underline{\rm End}_A(M)$
must be finite-dimensional, thus we get (3).
Then we copy the proof of \cite[Lemma 4.6]{Sc2} and use Lemma
\ref{stablemultiplicative} to show that each
path $\underline{f_{a_1}} \cdots \underline{f_{a_l}}$, which is 0
in $\underline{\rm End}_A(M)$, has to be of length 2.
As in \cite[Lemma 4.7]{Sc2} we get that there are no
commutativity relations.
Thus property (4) holds.
This finishes the proof.
\end{proof}

The following lemma is based on an idea of Ringel.

\begin{Lem}\label{uniserialimage}
Let $f_a: M(C_1) \to M(C_2)$ be a two-sided graph map with
$a = ((D_1,E_1,F_1),(D_2,E_2,F_2))$.
If $D_1E_1F_2$ and $D_2E_1F_1$ are strings, or if
$D_1E_1D_2^-$ and $F_2^-E_1F_1$ are strings, then
${\rm Ext}^1_A(M(C_2),M(C_1)) \not= 0$.
\end{Lem}

\begin{proof}
Without loss of generality assume $E_1 = E_2$, and set $E = E_1$.
Thus $D_1ED_2^-$ and $F_2^-EF_1$ cannot be strings.
If $D_1EF_2$ and $D_2EF_1$ are strings, we get a short exact sequence
\[
0 \longrightarrow M(D_1EF_1) \longrightarrow M(D_1EF_2) \oplus
M(D_2EF_1) \longrightarrow M(D_2EF_2) \longrightarrow 0
\]
of $A$-modules, see \cite{Sc2} for a precise construction.
Since $a$ is two-sided, it follows that
the middle term of this sequence is not isomorphic to the direct
sum of its end terms.
Thus the sequence does not split, which implies
${\rm Ext}_A^1(M(C_2),M(C_1)) \not= 0$.
\end{proof}

\begin{Cor}\label{coruniserialimage}
Let $f_a: M(C_1) \to M(C_2)$ be a two-sided graph map with
$a = ((D_1,E_1,F_1),(D_2,E_2,F_2))$.
If $E_1$ is not a direct or an inverse string, then
${\rm Ext}^1_A(M(C_2),M(C_1)) \not= 0$.
\end{Cor}

\begin{Prop}\label{factorizes}
Let $f_a: M(C_1) \to M(C_2)$ be a two-sided graph map, and assume
${\rm Ext}_A^1(M(C_2),M(C_1)) = 0$. Then $\underline{f_a} = 0$.
\end{Prop}

\begin{proof}
Let $a = ((D_1,E_1,F_1),(D_2,E_2,F_2))$.
Without loss of generality we can assume that $E_1 = E_2$.
Set $E = E_1$.
If $D_1EF_2$ and $D_2EF_1$ are both strings, then by
Lemma \ref{uniserialimage} we get ${\rm Ext}_A^1(M(C_2),M(C_1)) \not= 0$,
a contradiction.

Thus, assume that $D_2EF_1$ is not a string.
This implies that $E$ must be a direct string.

Let $F_{11}$ be a direct string of maximal length such that
$F_1 = F_{11}F_{12}$ for some string $F_{12}$.
Similarly, let $D_{22}$ be a maximal direct string such that
$D_2 = D_{21}D_{22}$ for some $D_{21}$.
Observe that $D_{22}EF_{11}$ is a string if and only if $D_2EF_1$
is a string.
Thus $D_{22}EF_{11}$ is not a string.
Since $D_{22}E$ and $EF_{11}$ are both strings, we get that
$|D_{22}| \geq 1$ and $|F_{11}| \geq 1$.
Thus $f_a: M(C_1) \to M(C_2)$ can be visualised as follows:

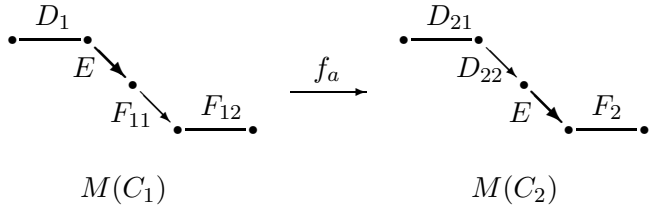
\unitlength1.0cm
\begin{picture}(12,3.3)

\put(0.5,2.5){\circle*{0.1}}
\put(0.6,2.5){\line(1,0){0.8}}\put(0.8,2.7){$D_1$}
\put(1.5,2.5){\circle*{0.1}}
{\thicklines
\put(1.6,2.4){\vector(1,-1){0.4}}\put(1.3,2){$E$}
}
\put(2.1,1.9){\circle*{0.1}}
\put(2.2,1.8){\vector(1,-1){0.4}}\put(1.8,1.4){$F_{11}$}
\put(2.7,1.3){\circle*{0.1}}
\put(2.8,1.3){\line(1,0){0.8}}\put(3,1.5){$F_{12}$}
\put(3.7,1.3){\circle*{0.1}}

\put(4.5,2){$f_a$}
\put(4.2,1.8){\vector(1,0){1.0}}

\put(5.7,2.5){\circle*{0.1}}
\put(5.8,2.5){\line(1,0){0.8}}\put(6,2.7){$D_{21}$}
\put(6.7,2.5){\circle*{0.1}}
\put(6.8,2.4){\vector(1,-1){0.4}}\put(6.4,2){$D_{22}$}
\put(7.3,1.9){\circle*{0.1}}
{\thicklines
\put(7.4,1.8){\vector(1,-1){0.4}}\put(7.1,1.4){$E$}
}
\put(7.9,1.3){\circle*{0.1}}
\put(8,1.3){\line(1,0){0.8}}\put(8.2,1.5){$F_2$}
\put(8.9,1.3){\circle*{0.1}}

\put(1.4,0.4){$M(C_1)$}
\put(6.6,0.4){$M(C_2)$}

\end{picture}

Now let $F_{111}$ be a maximal direct string such that $D_{22}EF_{111}$
is still a string.
Thus there exists a (direct) string $F_{112}$ with $|F_{112}| \geq 1$
such that $F_{11} = F_{111}F_{112}$.
Next, let $D_{212}$ be a maximal direct string such that
$D_{21} = D_{211}D_{212}^-$ for some string $D_{211}$.
Let $H$ be a maximal direct string such that $D_{212}H$ is
still a string.
Define $P = H^-D_{212}^-D_{22}EF_{111}$.
Since $|D_{22}| \geq 1$ and $|F_{112}| \geq 1$ we get a graph map
$f_b: M(C_1) \to M(P)$ with
\[
b = ((D_1,EF_{111},F_{112}F_{12}),(H^-D_{212}^-D_{22},EF_{111},
1_{(e(F_{111}),-\epsilon(F_{111}))})).
\]
The corresponding picture looks as follows:

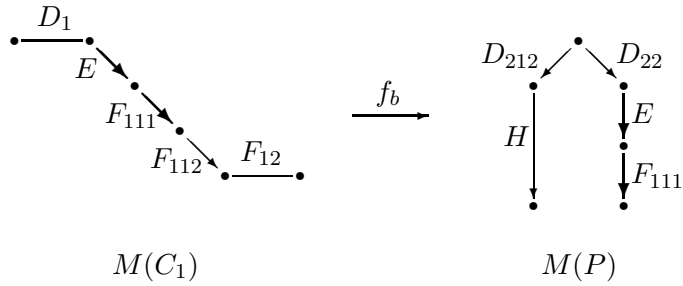
\unitlength1.0cm
\begin{picture}(12,4.3)

\put(0.5,3.5){\circle*{0.1}}
\put(0.6,3.5){\line(1,0){0.8}}\put(0.8,3.7){$D_1$}
\put(1.5,3.5){\circle*{0.1}}
{\thicklines
\put(1.6,3.4){\vector(1,-1){0.4}}\put(1.3,3){$E$}
}
\put(2.1,2.9){\circle*{0.1}}
{\thicklines
\put(2.2,2.8){\vector(1,-1){0.4}}\put(1.7,2.4){$F_{111}$}
}
\put(2.7,2.3){\circle*{0.1}}
\put(2.8,2.2){\vector(1,-1){0.4}}\put(2.3,1.8){$F_{112}$}
\put(3.3,1.7){\circle*{0.1}}
\put(3.4,1.7){\line(1,0){0.8}}\put(3.5,1.9){$F_{12}$}
\put(4.3,1.7){\circle*{0.1}}

\put(5.3,2.7){$f_b$}
\put(5,2.5){\vector(1,0){1.0}}

\put(8,3.5){\circle*{0.1}}
\put(8.1,3.4){\vector(1,-1){0.4}}\put(8.5,3.2){$D_{22}$}
\put(8.6,2.9){\circle*{0.1}}
{\thicklines
\put(8.6,2.8){\vector(0,-1){0.6}}\put(8.7,2.4){$E$}
}
\put(8.6,2.1){\circle*{0.1}}
{\thicklines
\put(8.6,2){\vector(0,-1){0.6}}\put(8.7,1.6){$F_{111}$}
}
\put(8.6,1.3){\circle*{0.1}}
\put(7.9,3.4){\vector(-1,-1){0.4}}\put(6.7,3.2){$D_{212}$}
\put(7.4,2.9){\circle*{0.1}}
\put(7.4,2.8){\vector(0,-1){1.4}}\put(7,2.1){$H$}
\put(7.4,1.3){\circle*{0.1}}

\put(1.8,0.4){$M(C_1)$}
\put(7.5,0.4){$M(P)$}

\end{picture}

Observe that $D_{211}$ is either of length 0 or ends with
an arrow.
Thus we get another graph map
$f_c: M(P) \to M(C_2)$ with
\[
c = ((H^-,D_{212}^-D_{22}E,F_{111}),(D_{211},D_{212}^-D_{22}E,F_2)),
\]
and the corresponding picture is the following:

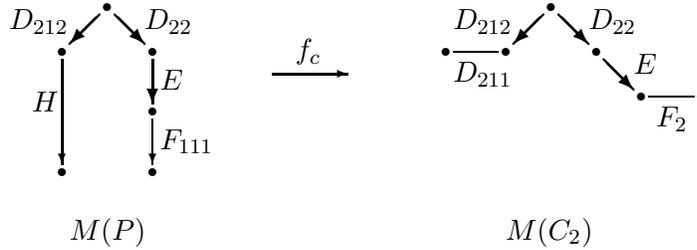
\unitlength1.0cm
\begin{picture}(12,3.9)

\put(1.6,3.5){\circle*{0.1}}
{\thicklines
\put(1.7,3.4){\vector(1,-1){0.4}}\put(2.1,3.2){$D_{22}$}
}
\put(2.2,2.9){\circle*{0.1}}
{\thicklines
\put(2.2,2.8){\vector(0,-1){0.6}}\put(2.3,2.4){$E$}
}
\put(2.2,2.1){\circle*{0.1}}
\put(2.2,2){\vector(0,-1){0.6}}\put(2.3,1.6){$F_{111}$}
\put(2.2,1.3){\circle*{0.1}}
{\thicklines
\put(1.5,3.4){\vector(-1,-1){0.4}}\put(0.3,3.2){$D_{212}$}
}
\put(1,2.9){\circle*{0.1}}
\put(1,2.8){\vector(0,-1){1.4}}\put(0.6,2.1){$H$}
\put(1,1.3){\circle*{0.1}}

\put(4.1,2.8){$f_c$}
\put(3.8,2.6){\vector(1,0){1.0}}

\put(7.5,3.5){\circle*{0.1}}
{\thicklines
\put(7.4,3.4){\vector(-1,-1){0.4}}\put(6.2,3.2){$D_{212}$}
}
\put(6.9,2.9){\circle*{0.1}}
\put(6.8,2.9){\line(-1,0){0.6}}\put(6.2,2.5){$D_{211}$}
\put(6.1,2.9){\circle*{0.1}}
{\thicklines
\put(7.6,3.4){\vector(1,-1){0.4}}\put(8,3.2){$D_{22}$}
}
\put(8.1,2.9){\circle*{0.1}}
{\thicklines
\put(8.2,2.8){\vector(1,-1){0.4}}\put(8.6,2.6){$E$}
}
\put(8.7,2.3){\circle*{0.1}}
\put(8.8,2.3){\line(1,0){0.6}}\put(8.9,1.9){$F_2$}
\put(9.5,2.3){\circle*{0.1}}

\put(1.1,0.4){$M(P)$}
\put(6.9,0.4){$M(C_2)$}

\end{picture}

Next, one checks easily that $f_a = f_bf_c$.
Thus $f_a$ factors through $M(P)$.
In case $M(P)$ is projective, we get $\underline{f_a} = 0$.

Otherwise, the indecomposable projective module $P_j$ with top $S_j$,
$j = s(D_{22})$,
must be a non-serial projective-injective module.
Assume we are in this case.
This implies that there exists a commutativity relation
$D_{212}H\alpha - D_{22}EF_{111}\beta$ in $\rho$ where $\alpha$ and
$\beta$ are arrows.
We know that $D_{212}H = \gamma K$ for some arrow $\gamma$ and
some (direct) string $K$, and
$D_{22} = \delta D_{22}'$ for some arrow $\delta$ and some (direct)
string $D_{22}'$.
Thus the radical ${\rm rad}(P_j)$ of $P_j$ is isomorphic to
$M(D_{22}'EF_{111}\beta\alpha^-K^-)$, and the
socle factor $P_j/{\rm soc}(P_j)$ is isomorphic to $M(P)$.
The next picture describes $P_j$, its radical and socle factor together
with the canonical inclusion and projection, respectively.

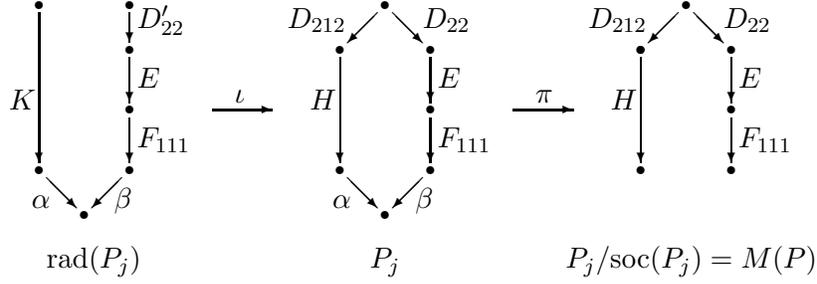
\unitlength1.0cm
\begin{picture}(12,4.3)

\put(1,3.9){\circle*{0.1}}
\put(1,3.8){\vector(0,-1){2}}\put(0.6,2.5){$K$}

\put(2.2,3.9){\circle*{0.1}}
\put(2.2,3.8){\vector(0,-1){0.4}}\put(2.3,3.6){$D_{22}'$}
\put(2.2,3.3){\circle*{0.1}}
\put(2.2,3.2){\vector(0,-1){0.6}}\put(2.3,2.8){$E$}
\put(2.2,2.5){\circle*{0.1}}
\put(2.2,2.4){\vector(0,-1){0.6}}\put(2.3,2){$F_{111}$}
\put(1,1.7){\circle*{0.1}}
\put(1.1,1.6){\vector(1,-1){0.4}}\put(0.9,1.2){$\alpha$}
\put(2.2,1.7){\circle*{0.1}}
\put(2.1,1.6){\vector(-1,-1){0.4}}\put(2,1.2){$\beta$}
\put(1.6,1.1){\circle*{0.1}}

\put(3.6,2.6){$\iota$}
\put(3.3,2.5){\vector(1,0){0.8}}

\put(5.6,3.9){\circle*{0.1}}
\put(5.5,3.8){\vector(-1,-1){0.4}}\put(4.3,3.6){$D_{212}$}
\put(5.7,3.8){\vector(1,-1){0.4}}\put(6.1,3.6){$D_{22}$}
\put(5,3.3){\circle*{0.1}}
\put(5,3.2){\vector(0,-1){1.4}}\put(4.6,2.5){$H$}

\put(6.2,3.3){\circle*{0.1}}
\put(6.2,3.2){\vector(0,-1){0.6}}\put(6.3,2.8){$E$}
\put(6.2,2.5){\circle*{0.1}}
\put(6.2,2.4){\vector(0,-1){0.6}}\put(6.3,2){$F_{111}$}
\put(5,1.7){\circle*{0.1}}
\put(5.1,1.6){\vector(1,-1){0.4}}\put(4.9,1.2){$\alpha$}
\put(6.2,1.7){\circle*{0.1}}
\put(6.1,1.6){\vector(-1,-1){0.4}}\put(6,1.2){$\beta$}
\put(5.6,1.1){\circle*{0.1}}

\put(7.6,2.6){$\pi$}
\put(7.3,2.5){\vector(1,0){0.8}}

\put(9.6,3.9){\circle*{0.1}}
\put(9.5,3.8){\vector(-1,-1){0.4}}\put(8.3,3.6){$D_{212}$}
\put(9.7,3.8){\vector(1,-1){0.4}}\put(10.1,3.6){$D_{22}$}
\put(9,3.3){\circle*{0.1}}
\put(9,3.3){\circle*{0.1}}
\put(9,3.2){\vector(0,-1){1.4}}\put(8.6,2.5){$H$}
\put(10.2,3.3){\circle*{0.1}}
\put(10.2,3.2){\vector(0,-1){0.6}}\put(10.3,2.8){$E$}
\put(10.2,2.5){\circle*{0.1}}
\put(10.2,2.4){\vector(0,-1){0.6}}\put(10.3,2){$F_{111}$}
\put(9,1.7){\circle*{0.1}}
\put(10.2,1.7){\circle*{0.1}}

\put(1.1,0.4){${\rm rad}(P_j)$}
\put(5.4,0.4){$P_j$}
\put(8,0.4){$P_j/{\rm soc}(P_j) = M(P)$}

\end{picture}

If $|F_{112}| > 1$, then it is easy to check that $f_a$ factors through
$P_j$.
Namely, we have $F_{112} = \beta F_{112}'$ for some (direct) string
$F_{112}'$ of length at least 1.
Then define $f_d: M(C_1) \to {\rm rad}(P_j)$
where
\[
d = ((D_1,EF_{111}\beta,F_{112}'F_{12}),
(D_{22}',EF_{111}\beta,\alpha^-K^-)).
\]
Thus $f_d$ looks as follows:

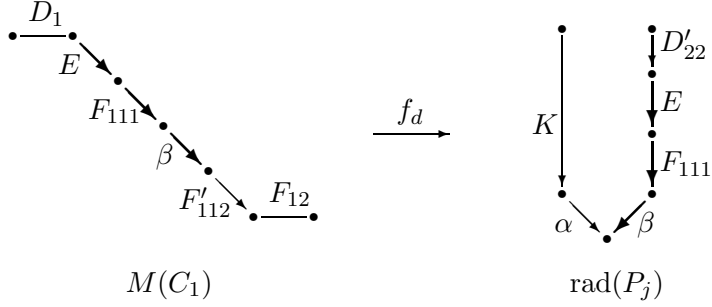
\unitlength1.0cm
\begin{picture}(12,4.5)

\put(0.2,3.8){\circle*{0.1}}
\put(0.3,3.8){\line(1,0){0.6}}\put(0.4,4){$D_1$}
\put(1,3.8){\circle*{0.1}}
{\thicklines
\put(1.1,3.7){\vector(1,-1){0.4}}\put(0.8,3.3){$E$}
}
\put(1.6,3.2){\circle*{0.1}}
{\thicklines
\put(1.7,3.1){\vector(1,-1){0.4}}\put(1.2,2.7){$F_{111}$}
}
\put(2.2,2.6){\circle*{0.1}}
{\thicklines
\put(2.3,2.5){\vector(1,-1){0.4}}\put(2.1,2.1){$\beta$}
}
\put(2.8,2){\circle*{0.1}}
\put(2.9,1.9){\vector(1,-1){0.4}}\put(2.4,1.5){$F_{112}'$}
\put(3.4,1.4){\circle*{0.1}}
\put(3.5,1.4){\line(1,0){0.6}}\put(3.6,1.6){$F_{12}$}
\put(4.2,1.4){\circle*{0.1}}

\put(5.3,2.7){$f_d$}
\put(5,2.5){\vector(1,0){1}}

\put(7.5,3.9){\circle*{0.1}}
\put(7.5,3.8){\vector(0,-1){2}}\put(7.1,2.5){$K$}
\put(8.7,3.9){\circle*{0.1}}
\put(8.7,3.8){\vector(0,-1){0.4}}\put(8.8,3.6){$D_{22}'$}
\put(8.7,3.3){\circle*{0.1}}
{\thicklines
\put(8.7,3.2){\vector(0,-1){0.6}}\put(8.8,2.8){$E$}
}
\put(8.7,2.5){\circle*{0.1}}
{\thicklines
\put(8.7,2.4){\vector(0,-1){0.6}}\put(8.8,2){$F_{111}$}
}
\put(7.5,1.7){\circle*{0.1}}
\put(7.6,1.6){\vector(1,-1){0.4}}\put(7.4,1.2){$\alpha$}
\put(8.7,1.7){\circle*{0.1}}
{\thicklines
\put(8.6,1.6){\vector(-1,-1){0.4}}\put(8.5,1.2){$\beta$}
}
\put(8.1,1.1){\circle*{0.1}}

\put(1.7,0.4){$M(C_1)$}
\put(7.6,0.4){${\rm rad}(P_j)$}

\end{picture}

Now one checks easily that
\[
f_a = f_d \iota \pi f_c.
\]
Thus $f_a$ factors through a projective module.

Next, assume $|F_{112}| = 1$.
This implies $F_{112} = \beta$.
Let $F_{121}$ be a maximal direct string such that
$F_{12} = F_{121}^-F_{122}$ for some string $F_{122}$.
Then there is a graph map $f_e: M(C_1) \to {\rm rad}(P_j)$ with
\[
e = ((D_1,EF_{111}\beta F_{121}^-,F_{122}),
(D_{22}',EF_{111}\beta F_{121}^-,R^-))
\]
where $R$ is a (direct) string such that
$RF_{121} = K\alpha$.
Thus $f_e$ looks as in the following picture:

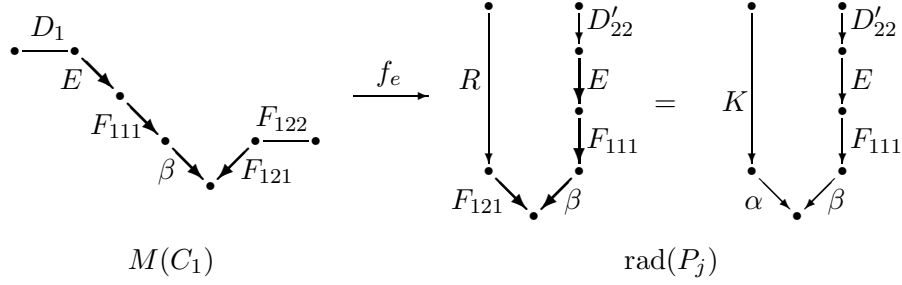
\unitlength1.0cm
\begin{picture}(12,4.3)

\put(0.2,3.3){\circle*{0.1}}
\put(0.3,3.3){\line(1,0){0.6}}\put(0.4,3.5){$D_1$}
\put(1,3.3){\circle*{0.1}}
{\thicklines
\put(1.1,3.2){\vector(1,-1){0.4}}\put(0.8,2.8){$E$}
}
\put(1.6,2.7){\circle*{0.1}}
{\thicklines
\put(1.7,2.6){\vector(1,-1){0.4}}\put(1.2,2.2){$F_{111}$}
}
\put(2.2,2.1){\circle*{0.1}}
{\thicklines
\put(2.3,2){\vector(1,-1){0.4}}\put(2.1,1.6){$\beta$}
}
\put(2.8,1.5){\circle*{0.1}}
{\thicklines
\put(3.3,2){\vector(-1,-1){0.4}}\put(3.2,1.6){$F_{121}$}
}
\put(3.4,2.1){\circle*{0.1}}
\put(3.5,2.1){\line(1,0){0.6}}\put(3.4,2.3){$F_{122}$}
\put(4.2,2.1){\circle*{0.1}}

\put(5,2.9){$f_e$}
\put(4.7,2.7){\vector(1,0){1}}

\put(6.5,3.9){\circle*{0.1}}
\put(6.5,3.8){\vector(0,-1){2}}\put(6.1,2.8){$R$}
\put(7.7,3.9){\circle*{0.1}}
\put(7.7,3.8){\vector(0,-1){0.4}}\put(7.8,3.6){$D_{22}'$}
\put(7.7,3.3){\circle*{0.1}}
{\thicklines
\put(7.7,3.2){\vector(0,-1){0.6}}\put(7.8,2.8){$E$}
}
\put(7.7,2.5){\circle*{0.1}}
{\thicklines
\put(7.7,2.4){\vector(0,-1){0.6}}\put(7.8,2){$F_{111}$}
}
\put(6.5,1.7){\circle*{0.1}}
{\thicklines
\put(6.6,1.6){\vector(1,-1){0.4}}\put(6,1.2){$F_{121}$}
}
\put(7.7,1.7){\circle*{0.1}}
{\thicklines
\put(7.6,1.6){\vector(-1,-1){0.4}}\put(7.5,1.2){$\beta$}
}
\put(7.1,1.1){\circle*{0.1}}

\put(8.7,2.5){$=$}

\put(10,3.9){\circle*{0.1}}
\put(10,3.8){\vector(0,-1){2}}\put(9.6,2.5){$K$}
\put(11.2,3.9){\circle*{0.1}}
\put(11.2,3.8){\vector(0,-1){0.4}}\put(11.3,3.6){$D_{22}'$}
\put(11.2,3.3){\circle*{0.1}}
\put(11.2,3.2){\vector(0,-1){0.6}}\put(11.3,2.8){$E$}
\put(11.2,2.5){\circle*{0.1}}
\put(11.2,2.4){\vector(0,-1){0.6}}\put(11.3,2){$F_{111}$}
\put(10,1.7){\circle*{0.1}}
\put(10.1,1.6){\vector(1,-1){0.4}}\put(9.9,1.2){$\alpha$}
\put(11.2,1.7){\circle*{0.1}}
\put(11.1,1.6){\vector(-1,-1){0.4}}\put(11,1.2){$\beta$}
\put(10.6,1.1){\circle*{0.1}}

\put(1.7,0.4){$M(C_1)$}
\put(8.3,0.4){${\rm rad}(P_j)$}

\end{picture}

If $|F_{121}| < |H\alpha|$, then one easily checks that
\[
f_a = f_e \iota \pi f_c.
\]
Thus, again $f_a$ factors through a projective.

It remains to consider the case $|F_{121}| \geq |H\alpha|$.
Note that this implies $|D_{212}| \geq 1$.
It follows that there is some
(direct) string $F_{1211}$ such that $F_{121} = F_{1211}F_{1212}$
and  $D_{212} = D_{2121}F_{1211}$
for some (direct) strings $F_{1212}$ and $D_{2121}$.
Note that $|F_{1212}| \geq 1$.
Thus we get a graph map $f_l: M(C_1) \to M(C_2)$ with
\[
l = ((D_1EF_{111}\beta F_{1212}^-,F_{1211}^-,F_{122}),
(D_{211},F_{1211}^-,D_{2121}^-D_{22}EF_2)).
\]
We get
\[
f_a + f_l = f_e \iota \pi f_c.
\]
Define $g = f_a + f_l$.
Clearly, we have $\underline{g} = 0$.
The following picture illustrates the situation:

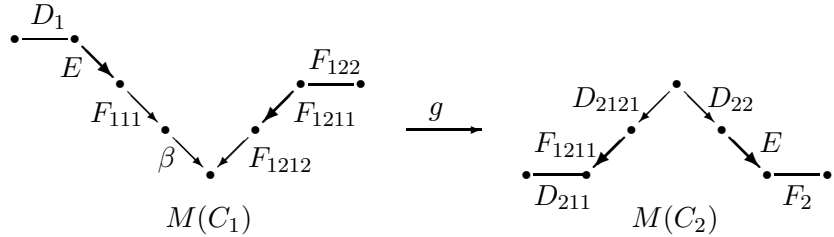
\unitlength1.0cm
\begin{picture}(12,3.2)

\put(0.2,2.5){\circle*{0.1}}
\put(0.3,2.5){\line(1,0){0.6}}\put(0.4,2.7){$D_1$}
\put(1,2.5){\circle*{0.1}}
{\thicklines
\put(1.1,2.4){\vector(1,-1){0.4}}}\put(0.8,2){$E$}
\put(1.6,1.9){\circle*{0.1}}
\put(1.7,1.8){\vector(1,-1){0.4}}\put(1.2,1.4){$F_{111}$}
\put(2.2,1.3){\circle*{0.1}}
\put(2.3,1.2){\vector(1,-1){0.4}}\put(2.1,0.8){$\beta$}
\put(2.8,0.7){\circle*{0.1}}

\put(3.3,1.2){\vector(-1,-1){0.4}}\put(3.3,0.8){$F_{1212}$}
\put(3.4,1.3){\circle*{0.1}}
{\thicklines
\put(3.9,1.8){\vector(-1,-1){0.4}}}\put(3.9,1.4){$F_{1211}$}
\put(4,1.9){\circle*{0.1}}
\put(4.1,1.9){\line(1,0){0.6}}\put(4.1,2.1){$F_{122}$}
\put(4.8,1.9){\circle*{0.1}}

\put(2.2,0){$M(C_1)$}

\put(5.7,1.5){$g$}
\put(5.4,1.3){\vector(1,0){1}}

\put(7,0.7){\circle*{0.1}}
\put(7.1,0.7){\line(1,0){0.7}}\put(7.1,0.3){$D_{211}$}
\put(7.8,0.7){\circle*{0.1}}
{\thicklines
\put(8.3,1.2){\vector(-1,-1){0.4}}}\put(7.1,1){$F_{1211}$}
\put(8.4,1.3){\circle*{0.1}}
\put(8.9,1.8){\vector(-1,-1){0.4}}\put(7.6,1.6){$D_{2121}$}
\put(9,1.9){\circle*{0.1}}

\put(9.1,1.8){\vector(1,-1){0.4}}\put(9.4,1.6){$D_{22}$}
\put(9.6,1.3){\circle*{0.1}}
{\thicklines
\put(9.7,1.2){\vector(1,-1){0.4}}}\put(10.1,1){$E$}
\put(10.2,0.7){\circle*{0.1}}
\put(10.3,0.7){\line(1,0){0.6}}\put(10.4,0.3){$F_2$}
\put(11,0.7){\circle*{0.1}}

\put(8.4,0){$M(C_2)$}

\end{picture}

\vspace{0.4cm}

Now we need some Auslander-Reiten theory:
If $|D_1| \geq 1$, then $f_a$ factors through $\tau_A^{-1} M(C_1)$.
This can be checked by using the construction of Auslander-Reiten
sequences for string algebras as explained in \cite{BR}.
Here we use also that $|F_{111}| \geq 1$.
Then we use the Auslander-Reiten formula
\[
{\rm Ext}_A^1(M(C_2),M(C_1)) \simeq
D \underline{\rm Hom}_A(\tau_A^{-1} M(C_1), M(C_2)).
\]
Thus, since ${\rm Ext}_A^1(M(C_2),M(C_1)) = 0$, we get
$\underline{f_a} = 0$.

Similarly, if $|F_{122}| \geq 1$, then $f_l$ factors through
$\tau_A^{-1} M(C_1)$.
Here we use that $|F_{1212}| \geq 1$.
Again, since ${\rm Ext}_A^1(M(C_2),M(C_1)) = 0$, we get
$\underline{f_l} = 0$.
But if $\underline{g} = \underline{f_l} = 0$, then
we get also $\underline{f_a} = 0$.

Thus, it remains to consider the case $|D_1| = |F_{122}| = 0$.
Then we can construct a short exact sequence
\[
0 \longrightarrow M(C_1) \longrightarrow
P_j \oplus M(D_{211}F_{1211}^-) \oplus M(EF_2)
\stackrel{h}{\longrightarrow} M(C_2)
\longrightarrow 0.
\]
Here
\[
h = [\pi f_c,f_m,f_n]^t
\]
where $f_m$ and $f_n$ are graph maps with
\[
m = ((1,D_{211}F_{1211}^-,1),(1,D_{211}F_{1211}^-,D_{2121}^-D_{22}EF_2))
\]
and
\[
n = ((1,EF_2,1),(D_{211}F_{1211}^-D_{2121}^-D_{22},EF_2,1)).
\]
For the sake of brevity,
in the definition of $m$ and $n$ we just wrote `1' for all
strings of length 0.
Next, one checks easily that $h$ is an epimorphism with kernel isomorphic
to $M(C_1)$.
This sequence does obviously not split,
thus ${\rm Ext}_A^1(M(C_2),M(C_1)) \not= 0$, a contradiction.
Altogether, we proved the following:
If $D_2EF_1$ is not a string, then the Proposition holds.

Thus, assume now that $D_2EF_1$ is a string.
If $D_1EF_2$ is also a string, then we apply Lemma \ref{uniserialimage}
and get that ${\rm Ext}_A^1(M(C_2),M(C_1)) \not= 0$, again a
contradiction.
Otherwise, if $D_1EF_2$ is not a string, then we proceed as before and get
that the Proposition holds.
This finishes the proof.
\end{proof}

\vspace{0.5cm}

\begin{proof}[Proof of Theorem \ref{main}]
Let $A$ be a special biserial algebra, and let $M$ be an $A$-module
with ${\rm Ext}_A^1(M,M) = 0$.
Thus $M$ does not contain direct summands which are isomorphic to
band modules.
Furthermore, since we are interested in stable endomorphism algebras,
we can assume that $M$ does not contain projective direct summands.
Also, since we consider algebras only up to Morita equivalence,
we assume that $M$ is a direct sum of pairwise different
indecomposable modules.
Thus $M$ is isomorphic to a direct sum $\bigoplus_{i=1}^t M(C_i)$
of pairwise different non-projective string modules $M(C_i)$.

The endomorphism algebra ${\rm End}_A(M)$ has a basis consisting of the
graph maps in ${\rm Hom}_A(M(C_i),M(C_j))$, $1 \leq i,j \leq t$.
By Proposition \ref{factorizes}
we know that any two-sided graph map $M(C_i) \to M(C_j)$ factorizes
through a projective module.
Hence we can apply Proposition \ref{newversion} and get that
$\underline{\rm End}_A(M)$ is a gentle algebra.
\end{proof}

\begin{proof}[Proof of Corollary \ref{maincor}]
If $B$ is a finite-dimensional $\field$-algebra, then
$\underline{\rm mod}(B)$ is a suspended category with suspension
functor $S = \Omega$ where $\Omega$ is the syzygy functor.
For $B$-modules $M$ and $N$ and $n \geq 1$ we have
\[
\underline{\rm Hom}_B(S^n(M),N) = {\rm Ext}_B^n(M,N)
\]
where the ${\rm Ext}_B^n(M,N)$ are the usual groups of module
extensions.
Now let $A$ be a finite-dimensional $\field$-algebra.
Then there is a
full and faithful embedding
\[
F: D^b(A) \longrightarrow \underline{\rm mod}(RA)
\]
of triangulated categories,
where $RA$ is the repetitive algebra of $A$.
For details we refer to \cite[Chapter 2]{H}
and \cite[Section 3]{K}.
For any two bounded complexes $M$ and $N$
of $A$-modules, we get
\begin{eqnarray*}
{\rm Hom}_{D^b(A)}(M[-1],N) &=&\underline{\rm Hom}_{RA}(F(M[-1]),F(N))\\
&=&\underline{\rm Hom}_{RA}(S(F(M)),F(N))\\
&=&{\rm Ext}^1_{RA}(F(M),F(N)),
\end{eqnarray*}
since $F$ is full, faithful and respects the triangulated
structure.

Now assume that $A$ is a finite-dimensional gentle algebra,
and let $T$ be a complex in $D^b(A)$
such that ${\rm Hom}_{D^b(A)}(T,T[1]) = 0$.
Using the above considerations,
we get ${\rm Ext}_{RA}^1(F(T),F(T)) = 0$.
Since $A$ is gentle, we know that $RA$ is special biserial.
Furthermore, ${\rm End}_{D^b(A)}(T)$ and
$\underline{\rm End}_{RA}(F(T))$ are isomorphic.
Thus we can apply Theorem \ref{main} and get that
${\rm End}_{D^b(A)}(T)$ is a gentle algebra.

If $T$ is a tilting complex, then by definition
${\rm Hom}_{D^b(A)}(T,T[i]) = 0$ for every $i \not= 0$.
But two algebras $A$ and $B$ are derived equivalent if and only if
there is a tilting complex $T$ in $D^b(A)$ such that $B$ is isomorphic to
${\rm End}_{D^b(A)}(T)$, see \cite{derbuch} or \cite{R}.
Thus any algebra which is derived equivalent to a
finite-dimensional gentle algebra is gentle again.

It is known that two derived equivalent finite-dimensional
algebras $A$ and $B$ have the
same number of isomorphism classes simple modules
\cite[Lemma 6.3.3]{derbuch}.
Since there are (up to Morita equivalence) only finitely many gentle
algebras with a given number of isomorphism classes of
simple modules, also the last statement in Corollary \ref{maincor} holds.
\end{proof}




\section{Examples and remarks}\label{examples}

\paragraph{\bf Example 1}
Let $A$ be the algebra $\field[x,y]/(x^2,y^2,xy)$ where $\field[x,y]$
is the polynomial algebra in two commuting variables.
Observe that $A$ is special biserial.
Now let $M$ be the string module $M(x^-yx^-y)$.
One easily checks that ${\rm Ext}_A^1(M,M) = 0$ and
${\rm End}_A(M) = \field[x_1, \cdots, x_6]/(x_ix_j, 1 \leq i,j \leq 6)$,
which is obviously not gentle and also far from being special biserial.
But we get
$\underline{\rm End}_A(M) = \field$.
Thus Theorem \ref{main} cannot be extended to a result on
endomorphism algebras, rather than stable endomorphism algebras.

\vspace{0.4cm}

\paragraph{\bf Example 2}
Now let $A = \field[x,y]/(x^3,y^3,xy)$ which is again special
biserial.
Let $M$ be the band module $M(x^-x^-yx^-yy,1,1)$, see \cite{BR}
for the construction of $M$.
Since $M$ is a band module we know that ${\rm Ext}_A^1(M,M) \not= 0$.
The endomorphism algebra ${\rm End}_A(M)$ is a factor of the free algebra
in 5 non-commuting variables,
and we get
$\underline{\rm End}_A(M) = \field$.
On the other hand, if we interpret $M$ as a module over the algebra
$B = \field[x,y]/(x^5,y^5,xy)$, then one easily checks that
${\rm End}_A(M) \simeq {\rm End}_B(M) \simeq \underline{\rm End}_B(M)$.

\vspace{0.4cm}

\paragraph{\bf Example 3}
The algebra $A = \field[x]/(x^2)$ is gentle, the simple module $S$
satisfies ${\rm Ext}_A^1(S,S) \not= 0$, and
${\rm End}_A(S) = \underline{\rm End}_A(S) = \field$ is gentle.

\vspace{0.4cm}

\paragraph{\bf Remark}
There are still several important open problems:
One needs criteria to decide when two given gentle algebras are
derived equivalent.

Also it is still not clear whether tilting-cotilting equivalence is the
same as derived equivalence.
This holds in the case of gentle one-cycle algebras, see \cite{AS}
and \cite{BGS}.





\begin{thebibliography}{80}

\bibitem{AH}
{\it I. Assem, D. Happel}, Generalized tilted algebras of type
$\mathbb{A}_n$, Comm. Algebra {\bf 9} (1981), no 20., 2101--2125.

\bibitem{ARS}
{\it M. Auslander, I. Reiten, S. Smal\o},
Representation theory of Artin algebras. Corrected reprint of the 1995 original,
Cambridge Studies in Advanced Mathematics, Vol. 36, Cambridge University
Press, Cambridge (1997), xiv+425pp.

\bibitem{AS}
{\it I. Assem, A. Skowro{\`n}ski}, Iterated tilted algebras of type
$\tilde{\mathbb{A}}_n$, Math. Z. {\bf 195} (1987), 269--290.

\bibitem{BGS}
{\it G. Bobi{\'n}ski, C. Gei{\ss}, A. Skowro{\'n}ski},
Classification of discrete derived categories, Preprint (2001),
1--35.

\bibitem{BR}
{\it M.C.R. Butler, C.M. Ringel}, Auslander-Reiten sequences with few middle
terms and applications to string algebras, Comm. Algebra {\bf 15} (1987),
no. 1-2, 145--179.

\bibitem{CB}
{\it W. Crawley-Boevey}, Maps between representations of zero-relation
algebras, J. Algebra {\bf 126} (1989), 259--263.

\bibitem{GP}
{\it I.M. Gelfand, V.A. Ponomarev}, Indecomposable representations
of the Lorentz group (Russian), Uspehi Mat. Nauk {\bf 23}
(1968), no. 2 (140), 3--60.

\bibitem{H}
{\it D. Happel}, Triangulated categories in the representation theory of
finite-dimensional algebras, London Mathematical Society Lecture Note
Series, Vol. 119, Cambridge University Press, Cambridge (1988), x+208pp.

\bibitem{K}
{\it B. Keller}, Chain complexes and stable categories, Manuscripta
Math. {\bf 67} (1990), no. 4, 379--417.

\bibitem{derbuch}
{\it S. K\"onig and A. Zimmermann}, 
Derived equivalences for group
rings, Lecture Notes in Mathematics, Vol. 1685, 
Springer-Verlag, Berlin (1998), x+246pp.

\bibitem{PS}
{\it Z. Pogorza{\l}y, A. Skowro{\'n}ski}, Selfinjective biserial standard
algebras, J. Algebra {\bf 138} (1991), 491--504.

\bibitem{R}
{\it J. Rickard}, Morita theory for derived categories,
J. London Math. Soc. (2) {\bf 39} (1989), no. 3, 436--456.

\bibitem{Ri1}
{\it C.M. Ringel}, Tame algebras and integral quadratic forms, In:
Lecture Notes in Mathematics, Vol. 1099, Springer-Verlag, Berlin
(1984), xiii+376pp.

\bibitem{Ri2}
{\it C.M. Ringel}, The repetitive algebra of a gentle algebra,
Bol. Soc. Mat. Mexicana (3) {\bf 3} (1997), no. 2, 235--253.

\bibitem{Sc1}
{\it J. Schr\"oer}, On the quiver with relations of a repetitive algebra,
Arch. Math. (Basel) {\bf 72} (1999), 426--432.

\bibitem{Sc2}
{\it J. Schr\"oer}, Modules without self-extensions over gentle algebras,
J. Algebra {\bf 216} (1999), 178--189.

\bibitem{Sc3}
{\it J. Schr\"oer}, On the infinite radical of a module category,
Proc. London Math. Soc. (3) {\bf 81} (2000), no. 3, 651--674.

\bibitem{Sc4}
{\it J. Schr\"oer}, On the Krull-Gabriel dimension of an algebra
Math. Z. {\bf 233} (2000), 287--303.

\bibitem{V}
{\it D. Vossieck}, The algebras with discrete derived category,
J. Algebra {\bf 243} (2001) 168-176.



\end{thebibliography}
\end{document}